\documentclass[preprint,12pt,3p]{elsarticle}

 \usepackage{graphicx}
 \usepackage{epstopdf}

\usepackage{amssymb}
\usepackage[final]{changes}
	\definechangesauthor[name={Per cusse}, color=orange]{per}
	\setremarkmarkup{(#2)}
\usepackage{amsmath}
\usepackage{xcolor}
\usepackage{caption}
 \usepackage{amsthm}

\usepackage{savesym}
\usepackage{amsmath}
\savesymbol{iint}
\usepackage{txfonts}
\restoresymbol{TXF}{iint}
\usepackage{booktabs}
\usepackage{hyperref}

\usepackage{subfig}
\usepackage{float}

\begin{document}

 \captionsetup[figure]{labelfont={bf},name={Fig.},labelsep=period}

\begin{frontmatter}

\title{Radial basis function collocation method for decoupled fractional Laplacian wave equations}

\author[Address1]{Yiran Xu}
\address[Address1]{ State Key Laboratory of Petroleum Resources and Prospecting,National Engineering Laboratory for Offshore Oil Exploration, China University of Petroleum,Beijing,102249, China}

\author[Address1]{Jingye Li\corref{cor1}}
\ead{ljy3605@sina.com}

\author[Address2]{Guofei Pang\corref{cor1}}
\address[Address2]{Algorithms Division,Beijing Computational Science Research Center, Beijing 100193, China}
\ead{pangguofei2008@126.com}

\author[Address1]{Zhikai Wang}


\author[Address1]{Xiaohong Chen}
\cortext[cor1]{The two corresponding authors contribute equally to the research.}

\begin{abstract}
Decoupled fractional Laplacian wave equation can describe the seismic wave propagation in attenuating media. Fourier pseudospectral implementations, which solve the equation in spatial frequency domain, are the only existing methods for solving the equation. For the earth media with curved boundaries, the pseudospectral methods could be less attractive to handle the irregular computational domains. In the paper, we propose a radial basis function collocation method that can easily tackle the irregular domain problems. Unlike the pseudospectral methods, the proposed method solves the equation in physical variable domain. The directional fractional Laplacian is chosen from varied definitions of fractional Laplacian. Particularly, the vector Gr\"unwald-Letnikov formula is employed to approximate fractional directional derivative of radial basis function. The convergence and stability of the method are numerically investigated by using the synthetic solution and the long-time simulations, respectively. The method's flexibility is studied by considering homogeneous and multi-layer media having regular and irregular geometric boundaries.

\end{abstract}


\end{frontmatter}

\section{Introduction}

Seismology, the scientific study of mechanical vibrations of the Earth, is used in mineral protecting and exploration for oil and natural gas, and in structural engineering to aid in the design of earthquake-resistant buildings. Conventional seismic modeling approaches either ignore the acoustic wave attenuation effects or use additional terms, which have a large number of parameters, to describe the attenuation effects. In modeling wave propagation in real media, amplitude loss and velocity dispersion effects need to be taken into account. Fractional derivative wave equations enjoy a fewer number of parameters and can describe the attenuation and the dispersion behaviors that exhibit frequency power-law dependency \cite{szabo1994time,chen2003modified,chen2004fractional}.

Fractional derivative modeling approaches are split into two categories: One category uses time-fractional derivatives that are introduced from the fractional derivative viscoelastic stress-strain relations \cite{holm2011causal,holm2014comparison}. The other \cite{chen2004fractional,treeby2010modeling,zhu2014modeling,Chen2017Fractional} takes advantage of the space-fractional derivative, namely fractional Laplacian, which is recovered by Fourier inverse transform from the frequency-domain equations under certain approximation assumption, say, the low-frequency assumption \cite{chen2004fractional}. Without a large memory of strain-stress history, the latter category compares favorably with the former one in terms of computational overhead when a long-time simulation is expected. There are still challenges for numerical simulation of fractional Laplacian wave equations.

First, varied definitions of fractional Laplacian have been introduced, and the guide to choose an optimal definition for a given application is absent. Notably, not all these definitions are equivalent on bounded domains.

Second,the effective numerical methods are needed to tackle the irregular domain wavefields. The Fourier pseudospectral methods \cite{carcione2010generalization,treeby2010modeling,zhu2014modeling,chen2016two} are shown to be rather efficient to simulate wavefield modeled by fractional Laplacian equations. However, the methods seem to be less attractive when irregular domains are considered, since the fast Fourier transform and its inversion are employed in the methods. The real boundary of seismic wavefield can be of arbitrary curvature, say, the hill or the valley on the earth surface.

Third, the discretization matrix of the fractional Laplacian is fully populated. Good preconditioning techniques deserve to be developed for solving the linear system having such a system matrix.

The present paper aims to solve the second problem aforementioned. The RBF collocation method \cite{Kansa1990Multiquadrics,Chen2014Recent,pang2015space-fractional,fornberg2015solving} can easily handle the discretization of high-dimensional, irregular computational domain, since the discretization depends only on the node-to-node distance. The method is mathematically simple and rather easy to program. Recently the method was for the first time applied to solve the advection-dispersion equation with the fractional Laplacian by the third author and his collaborators \cite{pang2015space-fractional}. In the present paper, an improved version of the method considered in \cite{pang2015space-fractional} is developed in order to solve a specific wave equation, called the decoupled fractional Laplacian wave equation \cite{zhu2014modeling}. In \cite{Anna2017}, the improved method has been used to solve two-dimensional fractional Poisson problems defined on bounded domains by the third author.

The paper is organized as follows. Sec.\ref{problem} introduces the decoupled fractional Laplacian wave equation as well as the definition of fractional Laplacian we are interested in. The RBF collocation method is elaborated in Sec.\ref{method}. Sec.\ref{result} shows the numerical simulations for homogeneous and multi-layer media of regular and irregular geometry boundaries. Concluding remarks are given in the last section.


\section{Problem}\label{problem}

\subsection{Nearly constant Q decoupled fractional Laplacian equation}

For homogeneous acoustic media, the nearly constant-$Q$ decoupled fractional Laplacian equation describing the stress field $\sigma(x,t)$ ($x\in\mathbb{R}^d$) is written by \cite{zhu2014modeling}
\begin{equation}
\frac   {1}  {c^2} \frac{\partial ^2 {\sigma}} { \partial ^2 {t}}=\eta(-\nabla^2)^{\gamma+1}\sigma+\tau\frac{\partial}{\partial{t}}(-\nabla^2)^{\gamma+\frac{1}{2}}\sigma
  \label{decoupled_fra}
\end{equation}
where $\eta=-c_0^{2\gamma}\omega_0^{-2\gamma}\cos(\pi\gamma)$, $\tau=-c_0^{2\gamma-1}\omega_0^{-2\gamma}\sin(\pi\gamma)$, and $c=c_0\cos({\pi\gamma}/{2})$ is the phase velocity depending on the velocity $c_0$ at the reference angular frequency $\omega_0$. The quality factor $Q$, which is related to $\gamma$ by $\gamma=\arctan(1/Q)/\pi$, is used to describe how fast the attenuation will be. The smaller the $Q$ is, the faster the wave attenuates. As $Q$ goes to the infinity, $\gamma$ will be zero, which indicates that Eq.(\ref{decoupled_fra}) will reduce to a standard wave equation (with no attenuation) due to $\tau = 0$ for $\gamma=0$. Oppositely, for $Q=0$, namely $\gamma=0.5$, the equation will be a diffusion equation, because $\eta$ is zero now. Note that for homogeneous media, $c_0$ and $Q$ are both independent of spatial variable and thus constants.

Under the condition $|\gamma\ln(\omega/\omega_0)|\ll 1$, the dispersive phase velocity $c$ and attenuation coefficient $\alpha$ arising from Eq.(\ref{decoupled_fra}) are given by
\begin{equation}
\begin{split}
  c & =c_0\left(\frac{\omega}{\omega_0}\right)^\gamma, \gamma\in (0,0.5),\\
  \alpha & =\tan\left(\frac{\pi\gamma}{2}\right)\frac{\omega}{c_p} \\
         & = \frac{\omega_0^{\gamma}\tan\left(\pi\gamma/2\right)}{c_0}\omega^{1-\gamma}.
\end{split}
\end{equation}
It is seen that the velocity and the attenuation are both frequency dependent for non-zero fractional order $\gamma$.

The reason why it is called the decoupled equation is that the equation can be decoupled into the attenuation-dominated equation
\begin{equation}
\frac{1}{c^2}\frac{\partial^2\sigma}{\partial^2 t^2}= \nabla^2\sigma+\tau\frac{\partial}{\partial t}(-\nabla^2)^{\gamma+1/2}\sigma, \gamma\in(0,0.5),
\end{equation}
and the dispersion-dominated equation
\begin{equation}
\frac{1}{c^2}\frac{\partial^2\sigma}{\partial^2 t^2}= \eta\frac{\partial}{\partial t}(-\nabla^2)^{\gamma+1}\sigma, \gamma\in(0,0.5).
\end{equation}
The decoupling effects are helpful for developing stable $Q$-compensated reverse time migration \cite{zhu2014q}.

\subsection{Fractional Laplacian}

There have been three typical definitions of the fractional Laplacian: integral or Riesz definition (see Eq.(25.56) of \cite{kilbas1993fractional}), directional definition (see \cite{meerschaert1999multidimensional,pang2015space-fractional} and Eq.(26.24) of \cite{kilbas1993fractional}), and spectral definition (\cite{ilic2005numerical,song2017computing}). Other definitions can be found in \cite{kwasnicki2017ten}. Ref.\cite{Anna2017} compares the aforementioned three typical definitions defined on bounded domains and concludes that the first two are equivalent but differ from the third one.

In the paper we restrict to the directional definition, given by \cite{pang2015space-fractional}
\begin{equation}\label{directional_def}
(-\nabla^2)^{\alpha/2}u(x)=C_{\alpha,d}\int_{||\boldsymbol{\theta}||=1}D_{\boldsymbol{\theta}}^{\alpha}u(x)d\boldsymbol{\theta}, \quad x,\boldsymbol{\theta} \in \mathbb{R}^d,
\end{equation}
where the scaling constant before the integral is \cite{pang2013gauss}
\begin{equation}
C_{\alpha,d}=\frac{\Gamma(\frac{1-\alpha}{2})\Gamma(\frac{d+\alpha}{2})}{2{\pi}^\frac{1+d}{2}},
\end{equation}
and the fractional directional derivative is given by
\begin{equation}
D_{\boldsymbol{\theta}}^{\alpha}(\cdot)=(\nabla \cdot \boldsymbol{\theta})^2 I_{\boldsymbol{\theta}}^{2-\alpha}(\cdot),
\end{equation}
where $\nabla$ is the gradient operator, and the fractional directional integral $I_{\boldsymbol{\theta}}^{\beta}(\cdot)$ is defined by (for $\beta\in(0,1)$)
\begin{equation}\label{fra_dir_int}
I_{\boldsymbol{\theta}}^{\beta}u(x)=\frac{1}{\Gamma(1-\beta)}\int_0^{+\infty}\varsigma^{-\beta}u(x-\varsigma\boldsymbol{\theta})d\varsigma,
\end{equation}
where $\Gamma(\cdot)$ is the Euler Gamma function. The directional definition is equivalent to the Riesz one and it is easy to be approximated by using RBF collocation method.

\subsection{Boundary and initial conditions}
Great care should be taken for setup of boundary conditions for the governing equation (\ref{decoupled_fra}). For pure Dirichlet problem, to guarantee uniqueness and existence of the solution, a nonlocal boundary condition should be given with the exterior condition and terminal observation time $T$ \cite{bucur2015some}:
\begin{equation}\label{zero-boundary}
\sigma(x,t)=0,\quad x\in \mathbb{R}^d\setminus \Omega, t\in [0,T].
\end{equation}
For simplicity, we here assume that the stress vanishes in the exterior. It should be noted that the above boundary condition could still be insufficient to ensure the well-posedness of the mathematical problem, since the highest order of the space derivative in our wave equation exceeds two because of the presence of the dispersion term $(-\nabla^2)^{\gamma+1}(\cdot)$ for $\gamma\in(0,0.5)$. An additional boundary condition could be needed. For example, the local-form normal condition
\begin{equation}\label{zero-normal}
\frac{\partial \sigma(x,t)}{\partial n}= 0, \quad x\in \partial \Omega,
\end{equation}
where $n$ is the outward unit normal, can be added. Other boundary conditions, say, the prescribed fractional normal derivative \cite{dipierro2014nonlocal}, is also possible. It should be noted that how to prove the well-posedness of the present wave equation given boundary conditions is still an open problem. Nevertheless, from the numerical experiments, we see that for zero boundary conditions, the normal boundary condition (\ref{zero-normal}) can be discarded. Therefore, in the rest of the paper, we will only adopt the zero boundary condition (\ref{zero-boundary}).

The initial conditions are given by
\begin{equation}
 \begin{split}
   \sigma(x,0) &= \sigma_0(x), \\
   \frac{\partial \sigma(x,0)}{\partial t} & = 0.
 \end{split}
\end{equation}

\section{Methodology}\label{method}
We first introduce the basic idea behind the RBF collocaiton method and then emphasize the computation of the fractional directional derivative of the RBFs using the vector Gr\"unwald-Letnikov formula. We restrict to two-dimensional problems and denote by $(x,y)$ the location of a point in two spatial dimensions. Note that in the preceding section the notation $x$ is used for representing the location of a point in the $d-$dimensional space.

\subsection{RBF collocation}
The stress field to be solved $\sigma(x,y,t)$ is approximated by the weighted sum of the RBFs $\phi(r)$:
\begin{equation}  \label{RBF-appr}
\begin{split}
 \sigma(x_i,y_i,t) & \approx  \sum_{j=1}^{M+N}{\lambda_j(t)\phi(r_{ij})},   \\
           r_{ij} & =  \sqrt{(x_i-x_j)^2+(y_i-y_j)^2},
 \end{split}
 \end{equation}
where $\{(x_j,y_j)\}$ are a group of source points located in the computational domain, and $M$ and $N$ are the numbers of points on the domain and boundary, respectively. The point set $\{(x_i,y_i)\}$ include collocation points, which coincide with the source points. The RBF $\phi(\cdot)$ only depends on the distance (or the relative location) between the collocation point $(x_i,y_i)$ and the source point $(x_j,y_j)$, namely $r_{ij}$, and therefore the collocation or source points are not constrained by any mesh or element. Additionally, $\lambda_j(t)$ is the time-dependent expansion coefficients to be evaluated.

Substituting RBF approximation (\ref{RBF-appr}) in our wave equation (\ref{decoupled_fra}) and in the boundary condition (\ref{zero-boundary}) and using the finite difference scheme for the temporal discretization, we lead to the following approximating scheme
\begin{equation} \label{semi-discrete}
\begin{split}
   \frac{1} {c_i^2} \sum_{j=1}^{M+N}\frac{\lambda_j^{n+1}-2\lambda_j^{n}+\lambda_j^{n-1}} {\Delta t^2}\phi_{ij} & = \eta\sum_{j=1}^{M+N}\lambda_j^{n}(-\nabla^2)^{\gamma+1}\phi_{ij}+\tau\sum_{j=1}^{M+N}\frac{\lambda_j^{n+1}-\lambda_j^n}{\Delta t}(-\nabla^2)^{\gamma+1/2}\phi_{ij}+f_i^n,  \\
    (x_i,y_i)\in \Omega,\quad i & = 1,2,\cdots,M, \\
   \sum_{j=1}^{M+N}\lambda_j\phi_{ij} & =0, \quad (x_i,y_i)\in \partial\Omega, \quad i=M+1,M+2,\cdots,M+N.
\end{split}
\end{equation}

The RBF expansion coefficient $\lambda_j^n$ means $\lambda_j((n-1)\Delta t)$ and the quantity $\phi_{ij}$ is defined by $\phi(r_{ij})$. $(-\nabla^2)^{\beta}\phi_{ij}$ represents $[(-\nabla^2)^{\beta}\phi](x_i,y_i)$ where $\beta$ can be $\gamma+1$ or $\gamma+1/2$. The first $M$ equations correspond to the approximation of the governing equation (\ref{decoupled_fra}) and the last $N$ equations approximate the zero-valued boundary condition (\ref{zero-boundary}). $\Delta t$ is the time step of the finite difference scheme. We here assume the velocity is spatially variable, i.e. $c_i=c(x_i,y_i)$, in order to use the model to simulate the layered media. An external force term $f_i^n=f(x_i,y_i,(n-1)\Delta t)$ is also considered. We have a trick here. The collocation (or source) points are placed simply on the boundary $\partial \Omega$ rather than the whole exterior $\mathbb{R}^d\setminus \Omega$. This trick avoids the choice of points outside $\Omega$. Without the trick, we need to place some points in the outward neighbourhood of boundary. To reduce the number of points, we remove the neighbourhood and find that the removing has little effect on the numerical results when zero-valued boundary condition is considered.

Rearranging the above approximating scheme, we have the matrix form ($n=2,3,...,L$ where $T=L\Delta t$)
\begin{equation}\label{time-iteration}
\left[
 \begin{array}{c}
 \boldsymbol{\Phi}_d-\tau \mathbf{C}^2\boldsymbol{\Phi}_{\gamma+1/2} \\
 \boldsymbol{\Phi}_b
 \end{array}
\right]\boldsymbol{\lambda}^{n+1} =
\left[
 \begin{array}{c}
 \eta \mathbf{C}^2\boldsymbol{\Phi}_{\gamma+1}-\tau \Delta t\mathbf{C}^2\boldsymbol{\Phi}_{\gamma+1/2} \\
 \boldsymbol{\Phi}_b
 \end{array}\right]\boldsymbol{\lambda}^n
 +\left[
 \begin{array}{c}
 -\boldsymbol{\Phi}_d \\
 \mathbf{0}
 \end{array}\right]\boldsymbol{\lambda}^{n-1}
 +\left[
 \begin{array}{c}
 \Delta t^2\mathbf{C}^2\mathbf{f}^n\\
 \mathbf{0}
 \end{array}
 \right].
\end{equation}
$\boldsymbol{\Phi}_d$ is a matrix formed by the RBFs with domain-type collocation points and all the source points, i.e. $[\phi_{ij}]$ with $i=1,2,\cdots,M; j=1,2,\cdots,M+N$. Similarly, $\boldsymbol{\Phi}_b$ is given by $[\phi_{ij}]$ with $i=M+1,M+2,\cdots,M+N; j=1,2,\cdots,M+N$. The $M\times M$ diagonal matrix $\mathbf{C}$ has the velocities evaluated at domain-type collocation points to be its diagonals, namely $[c_i]$ with $i=1,2,\cdots,M$. Corresponding to the $n$- th time layer, $\boldsymbol{\lambda}^n=[\lambda_j^n]$ and $\mathbf{f}^n=[f_i^n]$ are $(M+N)$- and $M$- dimensional column vectors, respectively. $\boldsymbol{\Phi}_{\gamma+1}$ and $\boldsymbol{\Phi}_{\gamma+1/2}$ are the $M\times (M+N)$ matrices formed by the fractional derivatives of the RBFs, and the next subsection will show how to compute these matrices.

The first expansion coefficient vector $\boldsymbol{\lambda}^1$ is obtained by solving the RBF interpolation problem
\begin{equation}
\left[
 \begin{array}{c}
  \boldsymbol{\Phi}_d \\
  \boldsymbol{\Phi}_b
  \end{array}
 \right]\boldsymbol{\lambda}^1 = \boldsymbol{\sigma}_0,
\end{equation}
where the vector formed by the initial values is denoted by $\boldsymbol{\sigma}_0=\sigma_0(x_i,y_i)$ for $i=1,2,\cdots,M+N$. The second expansion coefficient vector $\boldsymbol{\lambda}^2$ coincides with the first one due to the initial condition $\frac{\partial \sigma(x,y,0)}{\partial t}=0$, i.e.
\begin{equation}
 \boldsymbol{\lambda}^2=\boldsymbol{\lambda}^1.
\end{equation}
Through the iterations of the RBF expansion coefficients $\boldsymbol{\lambda}^n$ in (\ref{time-iteration}) starting from $\boldsymbol{\lambda}^1, \boldsymbol{\lambda}^2$, we can derive the coefficients for any given time $t_n=(n-1)\Delta t$. Finally, the use of RBF expansion formula (\ref{RBF-appr}) leads us to the approximate stress at any time-space coordinates.

\subsection{Fractional derivative of RBFs}\label{FD_RBF}

The computation of fractional Laplacian of RBFs is not usually an easy task, particularly when an unsuitable definition of fractional Laplacian is considered. For instance, if one considers the Riesz definition, which is a hyper-singular integral, even though the integrand is completely known, the hypersingularity prevents further approximation. Fortunately, the directional definition (\ref{directional_def}) enables a straightforward approximation of fractional Laplacian. The definition indicates that one can approximate the fractional Laplacian within two steps: Get the fractional directional derivative of the RBFs and then evaluate the integral with respect to differential direction $\boldsymbol{\theta}$.

Ref.\cite{pang2015space-fractional} approximates the fractional directional derivative of RBFs using Gauss-Jacobi quadrature under the assumption that the fractional directional derivative is defined in Caputo sense. Note that the Caputo fractional derivative is a regularized version of the Riemann-Liouville (RL) fractional derivative and does not completely preserve the nonlocality or the memory effect of the fractional derivative. For example, a RL fractional derivative of a constant is a power function, whereas a Caputo fractional derivative of a constant is zero, which is the same case as the integer-order derivative. Consideration of the RL fractional derivative could be preferable in mathematical modeling. To compute the RL derivative of RBFs, we adopt the vector Gr\"unwald scheme proposed in \cite{meerschaert2004vector}.

In two dimensional case, the scheme is written by
\begin{equation}\label{vGL-scheme}
D_{\theta}^{\beta}\phi(r_{ij})\approx h^{-\beta}\sum_{k=0}^{\left[\frac{d(x_i,y_i,\theta,\Omega)}{h}\right]}(-1)^k
\binom{\beta}{k}\phi\left(\sqrt{(x_i-kh\cos\theta-x_j)^2+(y_i-kh\sin\theta-y_j)^2}\right).
\end{equation}
The above scheme has the first order accuracy, i.e., $O(h)$ truncation error. The symbol $[a]$ gets the closest integer to the fraction $a$. $d(x_i,y_i,\theta,\Omega)$ is the distance from the collocation point $(x_i,y_i)$ to the boundary of the computing domain along the direction $(\cos\theta,\sin\theta)$. The second step is to evaluate the integration of the above directional derivative with respect to $\theta$. The integral can be accurately approximated by using trapezoidal rule since $D_{\theta}^{\beta}(\cdot)$ is periodic for $\theta$. We thus have the quadrature formula ($\theta_l=2\pi l/N_t$)
\begin{equation}
 \begin{split}
 (-\nabla^2)^{\beta}\phi(r_{ij}) & = C_{2\beta,2}\int_0^{2\pi}D_{\theta}^{\beta}\phi(r_{ij})d\theta \\
  & \approx \frac{2\pi C_{2\beta,2}}{N_t} \sum_{l=0}^{N_t-1}D_{\theta_l}^{\beta}\phi(r_{ij}),
 \end{split}
\end{equation}
where $N_t+1$ compound trapezoidal quadrature points are used. Letting the fractional order $\beta$ be $\gamma+1$ and $\gamma+1/2$, the index $i$ go from $1$ to $M$, and the index $j$ go from $1$ to $M+N$, we can compute all the elements of the matrices $\boldsymbol{\Phi}_{\gamma+1}$ and $\boldsymbol{\Phi}_{\gamma+1/2}$ aforementioned.

\section{Numerical results}\label{result}

In the section we first validate the RBF collocation method using synthetic solution. The convergence curve and the CPU time are shown. Second, for the square medium, we compare the solutions produced by RBF collocation and pseudo-spectral methods. The numerical stability of our method is also investigated via long-time simulation. Finally, the seismic wave RBF simulations in irregular, multi-layer media are highlighted, which illustrates the potentials of the proposed method for real seismic simulations.

As the frequently used RBF, the multi-quadratic function \cite{Kansa1990Multiquadrics} is taken as the RBF, namely $\phi(r)=\sqrt{r^2+p^2}$. The shape parameter $p$ is fixed to be the spatial step $\Delta x$. Although the RBF collocation method allows arbitrary, non-overlapping distribution of the collocation points, for convenience of programming, these points often coincide with the grid points of finite difference grid. Thus the spatial step $\Delta x$ is the very grid size. Additionally, in the following examples, the finite difference step for temporal discretization is set to be $\Delta t=10^{-7}s$.

\subsection{Validation of the RBF collocation method}

Consider the synthetic solution $\sigma_e=\exp(-t)x^3(1000-x)^3y^{3.6}(1000-y)^{3.6}$ on a square domain $\Omega=(0,1000)^2$. The external force term $f(x,y,t)$ can be directly computed by
\begin{equation}
f(x,y,t)=\frac{1}{c^2}\frac{\partial ^2 {\sigma_e(x,y,t)}} { \partial ^2{t}}-\eta(-\nabla^2)^{\gamma+1}\sigma_e(x,y,t)-\tau\frac{\partial}{\partial{t}}(-\nabla^2)^{\gamma+\frac{1}{2}}\sigma_e(x,y,t),
\end{equation}
with the velocity $c=3000 m/s$, the quality factor $Q=10$, the acoustic density $\rho_0=1 g/cm^3$, and the reference frequency $\omega_0=60 rad/s$. The computation of fractional Laplacian of a known function has been detailed in subsection \ref{FD_RBF}. The step size $h$ in the vector Gr\"unwald scheme (\ref{vGL-scheme}) is set to be $10^{-3}$ when computing the fractional Laplacian of RBFs. For computing the fractional Laplacian of $\sigma_e$, we halve the step size in order to make the external force $f$ as accurate as possible. Similarly, for computing the fractional Laplacian of RBFs, we use 21 trapezoidal quadrature points ($N_t=20$) while for the fractional Laplacian of $\sigma_e$, we use 41 quadrature points. Previous numerical experiments show that $h=1 m$ and $N_t=20$ can produce acceptable approximation accuracy. In the following examples, we keep these two parameters unchanged.

We test the RBF solution accuracy on the grid points of a $10\times 10$ regular grid on $\Omega$. The maximum absolute error and the average relative error are defined by $||\boldsymbol{\sigma}-\boldsymbol{\sigma}_e||_{\infty}$ and $||\boldsymbol{\sigma}-\boldsymbol{\sigma}_e||_2/||\boldsymbol{\sigma}_e||_2$, respectively. $\boldsymbol{\sigma}$ and $\boldsymbol{\sigma}_e$ are the vectors formed by the approximate solutions and the synthetic solutions evaluated at the test points.

\begin{figure}[H]
\centering
\includegraphics[width=0.7\textwidth]{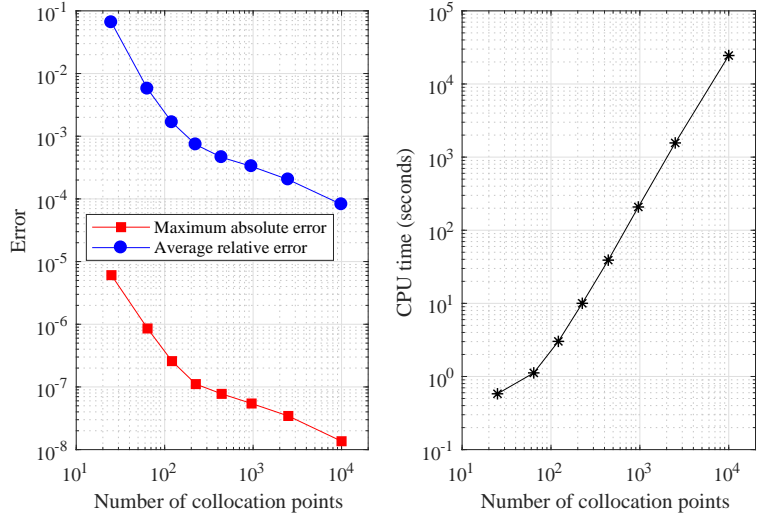}
\caption{\label{valid-RBF} Homogenous media with regular domain: Validation of the RBF collocation method (Error curve, left) and CPU time (right). All the solutions are evaluated at $t=10^{-5}s$ and the time step is $\Delta t=10^{-7}s$. To get the CPU time, a laptop with Intel core i7-5500U CPU (2.4G Hz) and 16 GB memory is used in computation. It is seen that the RBF method yields acceptable accuracy for a small number of collocation points, and the computational cost increases fast.}
\end{figure}

Fig.1(left) shows the variation of error against the number of collocation (or source) points, namely $M+N$. Fig.1(right) shows the total CPU time for computing the fractional Laplacian of RBFs, say,$\boldsymbol{\Phi}_{\gamma+1}$ and for directly solving the linear system (\ref{time-iteration}) once. We see that the RBF collocation method achieves acceptable accuracy for a small number of collocation points. But, the time complexity for generating the matrices whose entries are fractional Laplacian of RBFs is $O(N_tM(M+N))$ and the complexity for directly solving the linear system is $O((M+N)^3)$.

\subsection{RBF collocation versus Fourier pseudospectral method}

\begin{figure}[H]
\centering
\subfloat[Initial wavefield]{
\includegraphics[width=.45\textwidth]{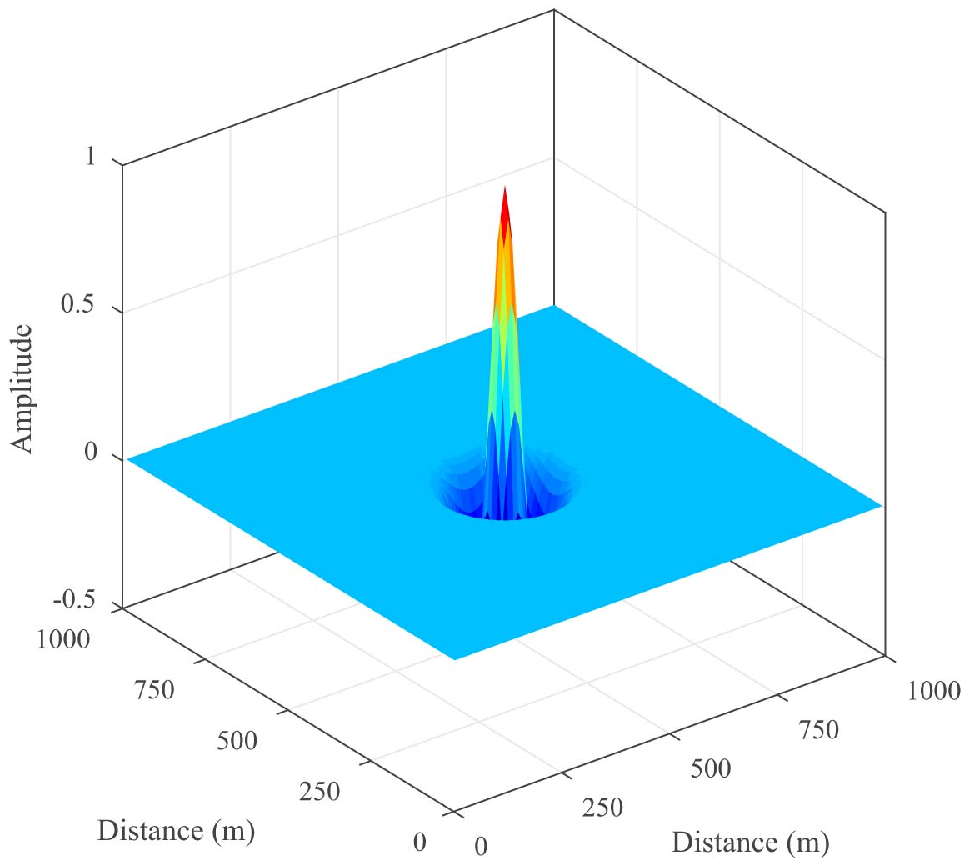}}\hfill
\centering
\subfloat[locations of 2601 collocation points (red dots for boundary points and blue dots for domain points)]{
\includegraphics[width=.4\textwidth]{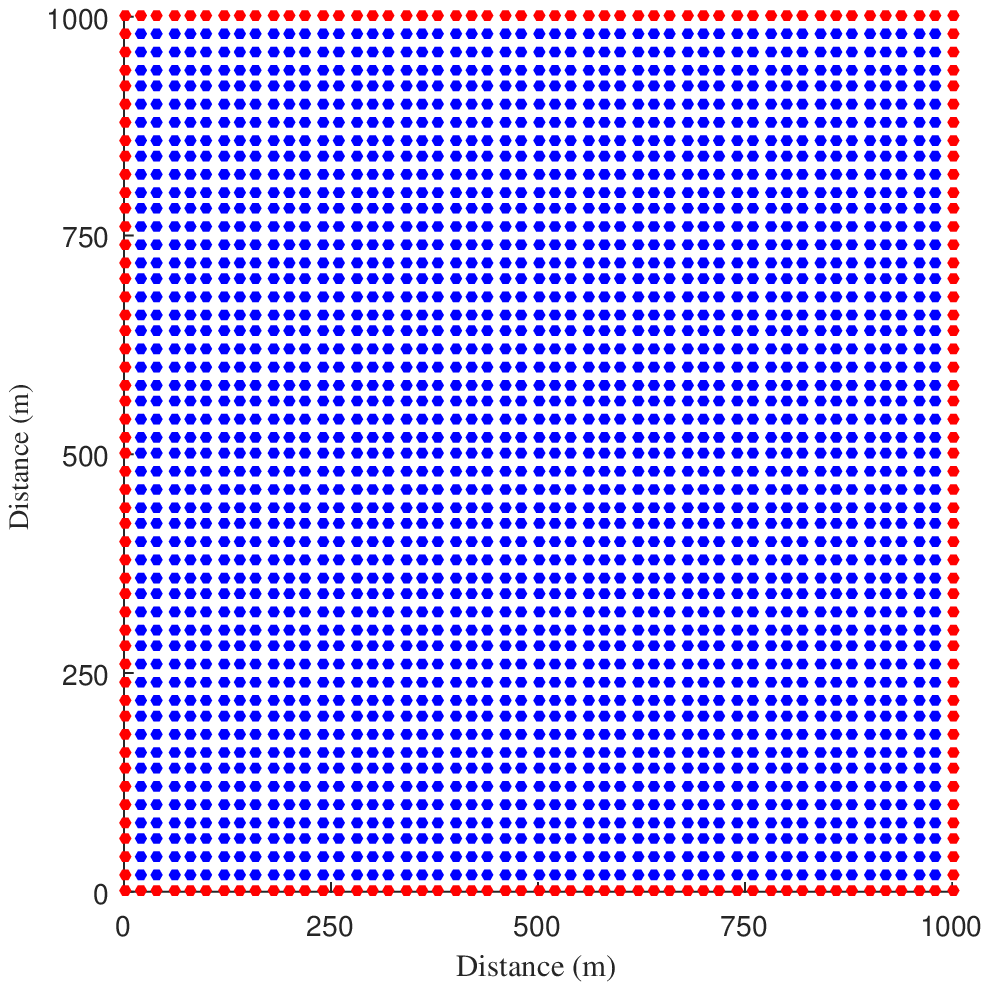}}
\caption{Homogenous media with regular domain: Initial wavefield (mimicking a point source) and the locations of collocation points used for computation.}
\label{initials}
\end{figure}

Consider a homogeneous medium of size $1000 m\times 1000 m$ with acoustic velocity $c_0=2000m/s$ and quality factor $Q=100$. Let the initial wavefield be
\begin{equation}
 \sigma(x,y,0)=[1-2(\pi f_0 r)^2]\exp[-(\pi f_0 r)^2]
 \end{equation}
 with $r^2=(x-x_s)^2+(y,y_s)^2$. The source location of the seismic wave is $(x_s,y_s)$. The parameter $f_0$ is set to be $f_0=5$. The plot of the initial wavefield is shown in Fig.\ref{initials}(a). In the rest of the paper, all the simulations use this initial wavefield, and $f_0=5$ and $f_0=12$ are used for homogeneous and multi-layer media, respectively. The external force is ignored, i.e. $f(x,y,t)\equiv 0$. The locations of collocation points are given in Fig.\ref{initials}(b).

\begin{figure}[H]
\centering
\subfloat[2D wavefield snapshot at $t=200$ ms for a homogenous model with acoustic reference velocity $c_0=2000m/s$ and quality factor $Q=100$: RBF method (left) and Pseudospectral method (right)]{
\includegraphics[width=.6\textwidth]{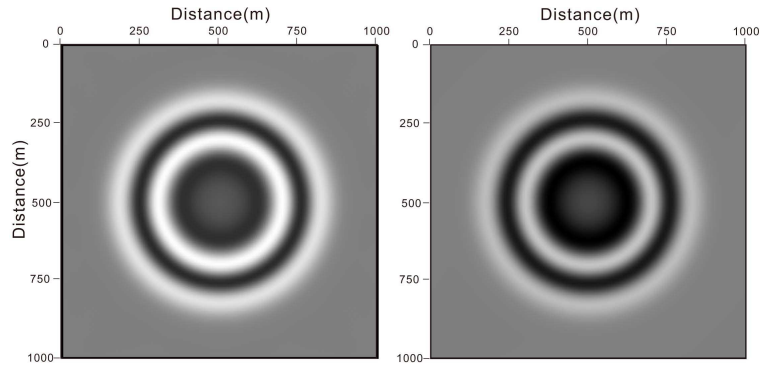}}\hfill
\subfloat[Amplitude along the line $y=500m$]{
\includegraphics[width=.7\textwidth]{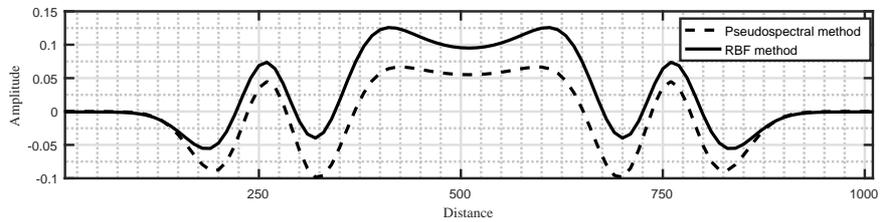}}
\caption{Homogenous media with regular domain: Comparison of RBF collocation and Fourier pseudospectral methods. Totally 10201 collocation points are used in the RBF method. Source location is $(500m, 500m)$. We see that the two methods produce the waves having the same phase but the different amplitudes.}
\label{wave-2d-RBF-PS-regular-homo}
\end{figure}

Fig.\ref{wave-2d-RBF-PS-regular-homo} compares the wavefields computed by RBF collocation and Fourier pseudospectral method. We see that the two methods produce the waves that have same phase but different amplitudes. The difference probably arises from the different boundary conditions considered in these two methods. The RBF collocation method considers the nonlocal zero-valued boundary condition, and the solutions on the exterior are zero; in contrast, the pseudospectral method assumes periodic boundary condition and the solutions on the exterior are simply the copies of the solution on the internal domain.

\begin{figure}[H]
\centering
\subfloat[Long-time wave propagation at $t=400 ms$]{
\includegraphics[width=.3\textwidth]{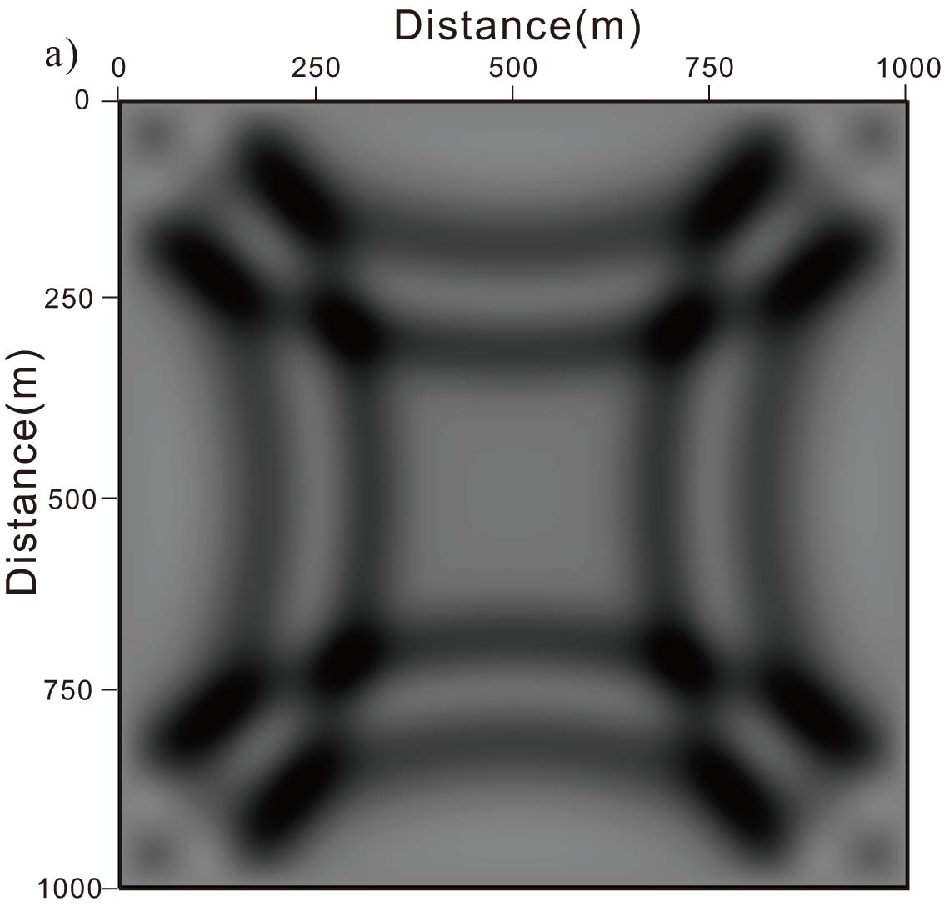}}\hfill
\subfloat[Long-time wave propagation at $t=4000 ms$]{
\includegraphics[width=.3\textwidth]{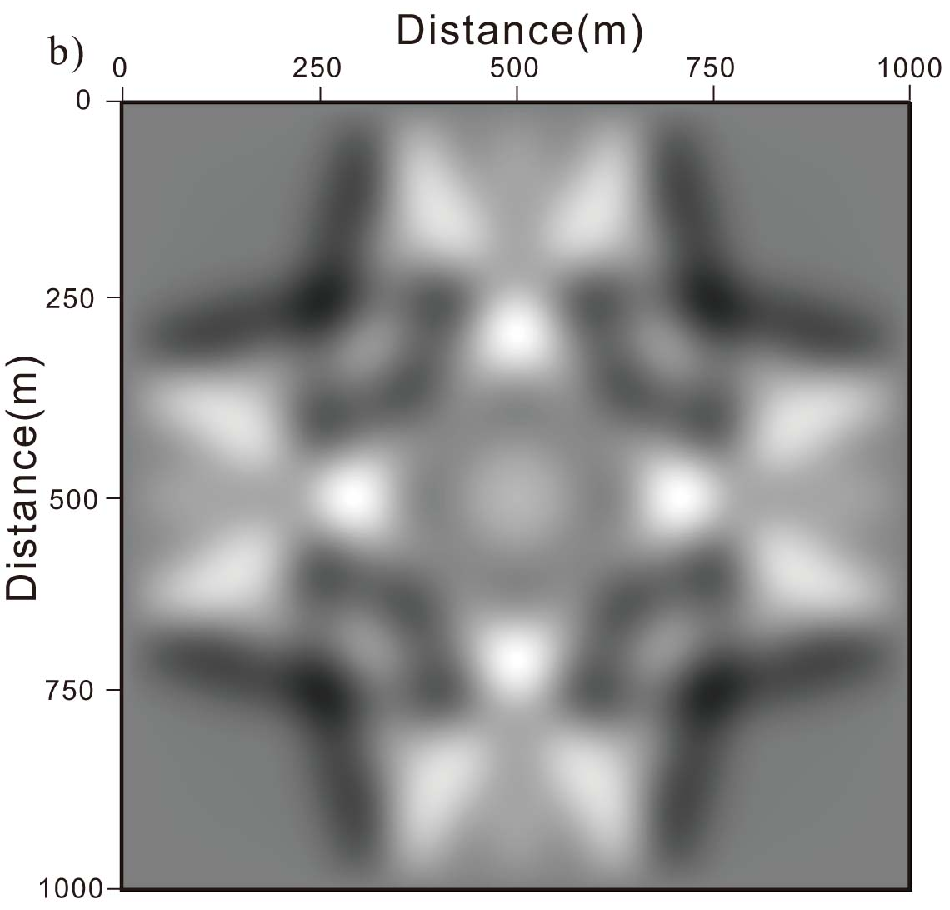}}\hfill
\subfloat[Simulated wavefields for varied quality factor $Q$ at $t=200 ms$]{
\includegraphics[width=.28\textwidth]{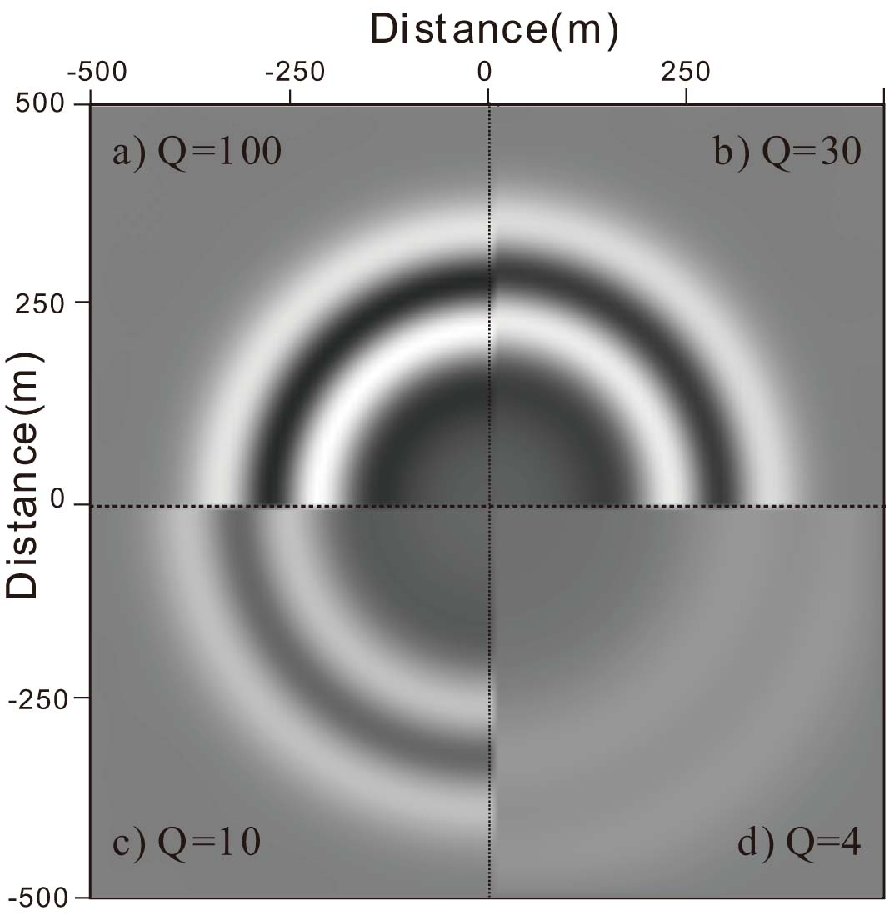}}
\caption{Homogenous media with regular domain: (a,b) Numerical stability of RBF method for long-time simulation; (c) Attenuation effects under different quality factors. Totally 10201 collocation points are used. We see that (1) The RBF method is stable and the corresponding solutions do not blow up for large $t$; (2) the smaller the quality factor $Q$ is, the larger the attenuation becomes and the faster the wave propagates. Source location is the center of the domain.}
\label{stability-Q}
\end{figure}

The numerical stability of the RBF collocation method for long-time simulations is studied and illustrated in Fig.\ref{stability-Q}(a,b). Additionally, the influence of the quality factor on the solutions is also considered, which can be observed from Fig.\ref{stability-Q}(c).

\subsection{Homogeneous media with irregular domain}
\begin{figure}[H]
\centering
\subfloat[Irregular domain I with 2091 collocation points]{
\includegraphics[width=.4\textwidth]{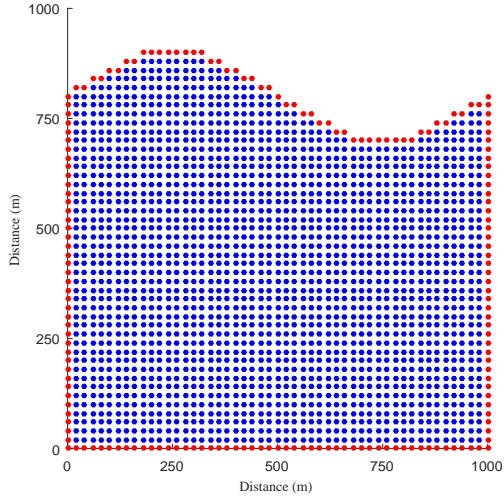}}\hfill
\subfloat[Irregular domain II with 2075 collocation points]{
\includegraphics[width=.4\textwidth]{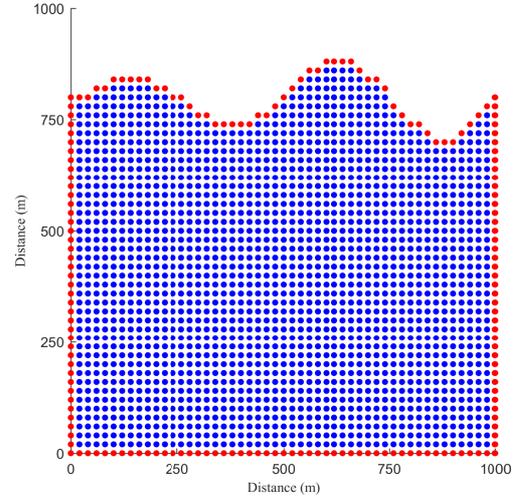}}\hfill
\subfloat[Simulated wavefield for domain I at time $t=220 ms$]{
\includegraphics[width=.4\textwidth]{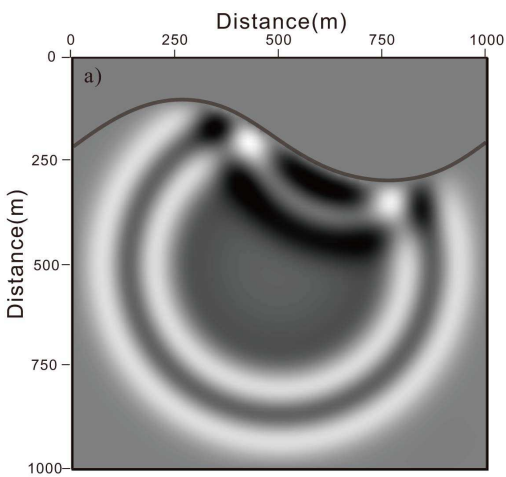}}\hfill
\subfloat[Simulated wavefield for domain II at time $t=220 ms$]{
\includegraphics[width=.4\textwidth]{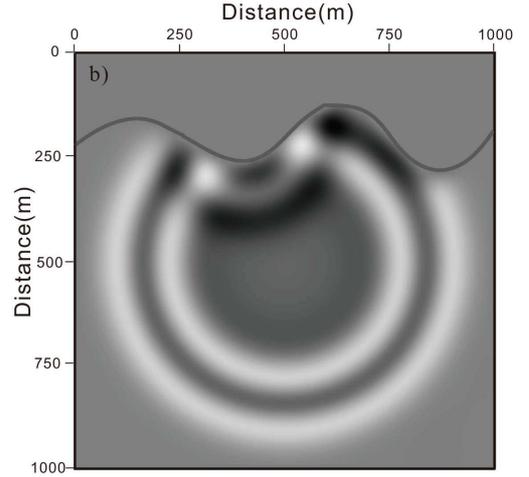}}
\caption{Homogenous media with irregular domain: RBF-simulated wavefields for two different irregular domains with the reference velocity $c_0=2000 m/s$ and the quality factor $Q=100$. Source is located on the center of the square domain. Totally 8181 and 8117 collocation points are used in RBF simulations for the first and the second domains, respectively. We see that the RBF method can handle the irregular domains of arbitrary geometric boundaries.}
\label{wave-homo-irregular}
\end{figure}

Fig. \ref{wave-homo-irregular} shows the wavefields simulated by using RBF method for two different irregular domains. These results are only used to illustrate the utility of the RBF collocation method for solving irregular domain problem. It should be noted that we make the first attempt to solve the decoupled fractional Laplacian equation on irregular domains.

\subsection{Multi-layer media}
\begin{figure}[H]
\centering
\includegraphics[width=.6\textwidth]{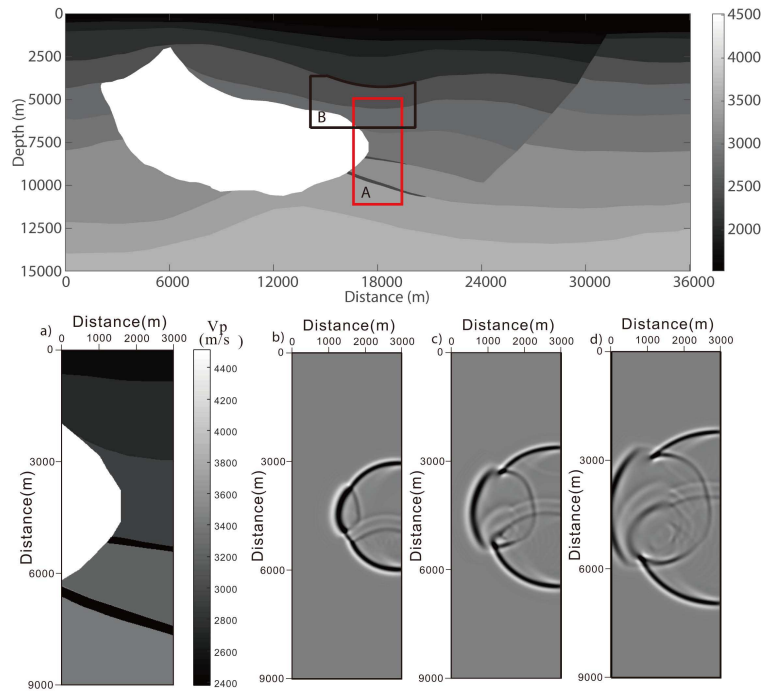}
\caption{Multi-layer media with regular domain: Synthetic phase velocity field of the original HESS VTI model (top subplot); Regions A and B are the regular and irregular subdomains we are interested in, respectively. In the second row, the leftmost subplot is the velocity distribution and the remaining ones are RBF-simulated wavefields at time $t=400,500,600 ms$. The quality factor is fixed to be a constant, i.e., $Q=100$. Source is located on the center of the right side. Totally $30401$ collocation points are used. We see that the RBF method can also tackle the one-direction heterogenous medium.}
\label{wave-layered-regular}
\end{figure}

We consider the media with acoustic velocity varying layer by layer. The synthetic velocity data are extracted from the synthetic dataset of the HESS VTI model \cite{liu2009decoupled}. The geometric size of the medium in the original HESS model is quite large, and to reduce the computational cost of solving the fractional Laplacian wave equation, we select a small part of medium. Fig.\ref{wave-layered-regular} (top) shows a $36000m\times 15000m$ medium considered in HESS model and what we consider are two subregions: Region A with regular boundary and Region B with irregular boundary.

The second row of the subfigures in Fig.\ref{wave-layered-regular} displays, respectively, the synthetic velocity field of Region A, and the simulated wavefields at $t=400, 500, 600 ms$. The initial wave source is located at the center of right side. The reflection and refraction of the wave on the interfaces of different layers can be clearly seen.

\begin{figure}[H]
\centering
\includegraphics[width=.9\textwidth]{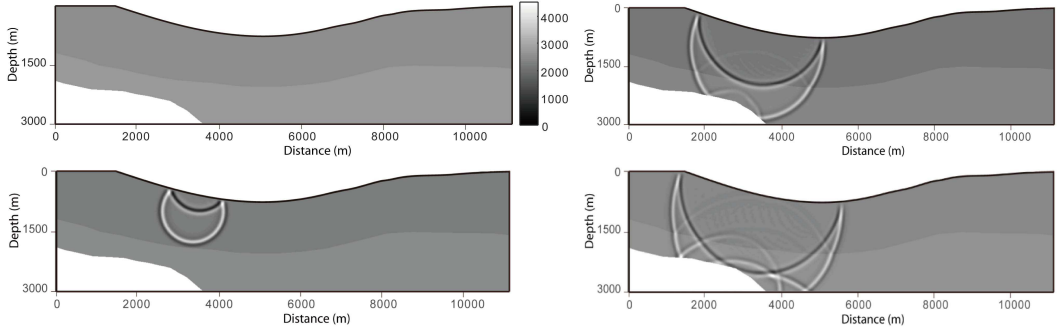}
\caption{Multi-layer media with irregular domain: Synthetic phase velocity field (upper left), RBF-simulated wavefield at $t=400 ms$ (lower left), $t=600 ms$ (upper right), and $t=800 ms$ (lower right). The quality factor is fixed to be $Q=100$. Source is located on the curved boundary. Totally $34165$ collocation points are used. It can be seen that RBF method can handle both the parameter (say, velocity) heterogeneity and the curved boundaries.}
\label{wave-layered-irregular}
\end{figure}

Fig.\ref{wave-layered-irregular} shows the simulations on Region B. The source is located on the curved boundary. The left top subfigure shows the synthetic velocity distribution, and the rest of the subfigures give the simulated wavefields at $t=400, 600, 800 ms$.

It can be seen that the RBF collocation method can be employed to simulate the wave in the media with spatially dependent velocity.

\section{Concluding remarks}

The paper proposes the radial basis collocation method for solving decoupled fractional Laplacian wave equation. The proposed method is easy to program and can simulate the wave propagation in media with irregular geometry boundaries. The synthetic solution was used to test the convergence of the method, and the synthetic spatially dependent velocity data on the HESS VTI model were used to show the flexibility of the method. The method could have certain potentials in forward seismic modeling.

High computational cost hinders the use of the method for large-scale simulations. It is desirable to develop efficient preconditioned iteration solver for the linear system generated by the method.

Variable fractional-order $\gamma(x,y)$ is another obstacle before increasing the popularity of the space-fractional derivative seismic forward modeling. The variable-order fractional Laplacian equations need to be taken into account; in fact, the quality factor of the real media is generally spatially dependent, and the fractional order $\gamma$ is related to the quality factor $Q$ by $\gamma=\arctan(1/Q)/\pi$. This requires the fractional order to also be spatially dependent. Unfortunately, the theoretical analysis of this type of equations lags far behind the requirements of using the equations in engineering world. It deserves to mention that in the present paper, when simulating the waves in multi-layer media, we keep the quality factor $Q$ to be a constant despite the spatially varied velocity, in order to avoid the variable-order modeling. We will exclusively discuss the RBF collocation method for variable-order equations in the coming work.

Boundary conditions, such as perfect match layer \cite{berenger1996perfectly}, need also to be investigated in the future.

\section{Acknowledgements}
This work was financially supported by National Natural Science Foundation of China (41774129, 41774131), National Science and Technology Major Project (2016ZX05024001-004), Science and Technology Project of CNPC (g2016A-3303), National Basic Research Program of China (973 Project No. 2010CB832702). The work of the third author was supported by the National Natural Science Foundation of China (11701025).

\bibliographystyle{elsarticle-num}

\bibliography{reference}

\end{document}